\documentclass[11pt]{amsart}
\usepackage{amsmath}
\usepackage{epic,eepic}
\advance\textwidth 3cm
\advance\oddsidemargin -1.5cm
\advance\evensidemargin -1.5 cm
\advance\textheight 2cm
\advance\topmargin -1cm
\title{Note on $\Pi^0_{n+1}$-LEM, $\Sigma^0_{n+1}$-LEM and $\Sigma^0_{n+1}$-DNE}
\author{\sc{Joan R. Moschovakis} \\ M.$\Pi.\Lambda$.A., \sc{Athens} \\ \sc{Emerita Math Prof., Occidental College}}
\begin{document}

\def\bibauthoryear{\begingroup
\def\thebibliography##1{\section*{References}%
    \small\list{}{\settowidth\labelwidth{}\leftmargin\parindent
    \itemindent=-\parindent
    \labelsep=\z@
    \usecounter{enumi}}%
    \def\newblock{\hskip .11em plus .33em minus -.07em}%
    \sloppy
    \sfcode`\.=1000\relax}%
    \def\@cite##1{##1}%
    \def\@lbibitem[##1]##2{\item[]\if@filesw
      {\def\protect####1{\string ####1\space}\immediate
    \write\@auxout{\string\bibcite{##2}{##1}}}\fi\ignorespaces}%
\endgroup}

\begin{abstract}
In \cite{ABHK} Akama, Berardi, Hayashi and Kohlenbach used a monotone modified realizability
interpretation to establish
the relative independence of $\Sigma^0_{n+1}$-DNE from $\Pi^0_{n+1}$-LEM over {\bf HA}, and hence the independence of
$\Sigma^0_{n+1}$-LEM from $\Pi^0_{n+1}$-LEM over {\bf HA}, for all $n \geq 0$.  We show that the same relative
independence results hold for these arithmetical principles over Kleene and Vesley's system {\bf FIM} of intuitionistic
analysis \cite{kleenevesley},
which extends {\bf HA} and is consistent with {\bf PA} but not with classical analysis.\footnote{Not even with
$\forall \alpha [\forall x (\alpha(x) = 0) \vee (\neg \forall x \alpha(x) = 0)]$.
In contrast, the extension of Markov's Principle ($\Sigma^0_1$-DNE) to the two-sorted
language {\it is} consistent with {\bf FIM}.  In {\bf FIM} + MP (but not in {\bf FIM}) it is possible to prove that
the constructive arithmetical hierarchy is proper; cf. \cite{myhierarchiespaper}, which also shows that {\bf FIM} is
not conservative over {\bf HA} with respect to arithmetical formulas.}
The double negations of the closures of $\Sigma^0_{n+1}$-LEM, $\Sigma^0_{n+1}$-DNE  and $\Pi^0_{n+1}$-LEM are also
considered, and shown to behave differently with respect to {\bf HA} and {\bf FIM}.  Various elementary questions
remain to be answered.
\end{abstract}

\thanks{I am grateful to Ulrich Kohlenbach for pointing me to \cite{ABHK}, and to the organizers of the
2005 Oberwolfach conference on Proof Theory and Constructive Mathematics for a terrific mathematical experience. \\
$\hfil^1${\bf FIM} is not even consistent with
$\forall \alpha [\forall x (\alpha(x) = 0) \vee (\neg \forall x \alpha(x) = 0)]$.
In contrast, the extension of Markov's Principle ($\Sigma^0_1$-DNE) to the two-sorted
language {\it is} consistent with {\bf FIM}.  In {\bf FIM} + MP (but not in {\bf FIM}) it is possible to prove that
the constructive arithmetical hierarchy is proper; cf. \cite{myhierarchiespaper}, which also shows that {\bf FIM} is
not conservative over {\bf HA} with respect to arithmetical formulas.}

\maketitle

{\em Definitions of the Arithmetical Principles}.
Unless otherwise noted, ``LEM'' (Law of Excluded Middle), ``DNE'' (Double Negation Elimination), and ``LLPO''
(Lesser Limited Principle of Omniscience) denote the (universal closures of the) purely arithmetical schemas,
without function variables.  If $\Phi$ is $\Sigma^0_n$ or $\Pi^0_n$ for some $n \geq 1$ then

(i) $~\Phi$-LEM is  $~A \vee \neg A~$  where $A \in \Phi$.

(ii) $~\Phi$-DNE is  $~\neg \neg A \rightarrow A~$ where $A \in \Phi$.

(iii) $~\Phi$-LLPO is $~\neg (A \wedge B) \rightarrow (C \vee D)~$, where $A,B \in \Phi$ and $C,D$ are the
duals of $A,B$ respectively.

(iv) $~\Delta^0_n$-LEM is $~(A \leftrightarrow B) \rightarrow (B \vee \neg B)$ where
$A \in \Pi^0_n$ and $B \in \Sigma^0_n$.

The precise statement of $\Delta^0_n$-LEM is important, since
$\Sigma^0_{n+1}$-DNE is equivalent over {\bf HA} + $\Sigma^0_n$-LEM to the schema
$~(\neg A \leftrightarrow B) \rightarrow (A \vee \neg A)~$ where $A, B \in \Sigma^0_{n+1}$.  Kleene used this principle
for $n=0$ to prove that every $\Delta^0_1$ relation is recursive.  The corresponding observation for $n \geq 0$
is the Kleene-Post-Mostowski Theorem.

\section{Some Results of Akama, Berardi, Hayashi and Kohlenbach Extended to {\bf FIM}}

{\em Lemma 1}.  The following are equivalent, for any theory {\bf T} $\supseteq$ {\bf HA}:
\begin{enumerate}
\item[(i)]  {\bf T} + $\Pi^0_1$-LEM  proves $\Sigma^0_1$-LEM.
\item[(ii)]  {\bf T} + $\Pi^0_1$-LEM  proves Markov's Principle $\Sigma^0_1$-DNE.
\end{enumerate}

{\it Proof}.  $(i) \Rightarrow (ii)$ holds because decidable predicates are stable under double negation.
$(ii) \Rightarrow (i)$ holds because
$$[\forall x \neg R(x) \vee \neg \forall x \neg R(x)] ~\&~ [\neg \neg \exists x R(x) \rightarrow \exists x R(x)]
\rightarrow [\exists x R(x) \vee \neg \exists x R(x)]\enspace.$$

Now let $T(e,x,y)$ be a quantifier-free formula numeralwise expressing in {\bf HA} (hence also in {\bf FIM}) the
Kleene T-predicate, and let $z \leq U(y)$ be a quantifier-free formula numeralwise expressing in {\bf HA} (hence
also in {\bf FIM}) the relation ``z $\leq$ U(y)'' where U(y) is the value computed by the computation with g\"odel
number y, or the g\"odel number of y if y is not the g\"odel number of a computation.  With Kleene's coding
{\bf HA} proves $\forall e \forall x \forall y [T(e,x,y) \rightarrow \forall z (z \leq U(y) \rightarrow
\neg T(e,x,z))]$, and we will use this property to prove the next lemma.

\vskip 0.1cm

{\it Lemma 2}.  {\bf HA} (hence also {\bf FIM}) proves
$$\forall f \neg \forall x \exists y [T(f,x,y) \wedge [\forall z_{z \leq U(y)}\neg T(x,x,z) \rightarrow
\forall y \neg T(x,x,y)]]\enspace.$$

{\it Proof}.  Assume for contradiction
$$\forall x \exists y [T(f,x,y) \wedge [\forall z_{z \leq U(y)}\neg T(x,x,z) \rightarrow \forall y \neg T(x,x,y)]]\enspace.$$
After $\forall$-elimination assume for $\exists y$-elimination:
$$T(f,f,y) \wedge [\forall z_{z \leq U(y)}\neg T(f,f,z) \rightarrow \forall y \neg T(f,f,y)]\enspace,$$
from which $T(f,f,y) \wedge \forall y \neg T(f,f,y)$ follows by the remark on coding.

\vskip 0.1cm

{\bf FIM} satisfies the ``independence of (stable) premise'' rule IPR:
$$(\ast) \; \mbox{If } \vdash_{\bf FIM} (\neg A \rightarrow \exists x B(x))  \mbox{ then }
\vdash_{\bf FIM} \exists x (\neg A \rightarrow B(x))$$
where $x$ is not free in $A$.  The beautiful proof by Visser that {\bf HA} is closed under IPR (cf. p. 138 of \cite{TvD})
works also for {\bf FIM}.  If one uses the monotone form ($^*$27.13 in \cite{kleenevesley}) of the bar induction schema,
it is straightforward to show that {\bf FIM} proves the Friedman translation of each of its mathematical axioms, and the
logical rules and axioms behave as usual.

\vskip 0.1cm

{\it Lemma 3}.  {\bf FIM} + $\Pi^0_1$-LEM  does not prove $\Sigma^0_1$-LEM.

\vskip 0.1cm

{\it Proof}.  We use without much comment the fact that quantifier-free formulas are decidable and stable in {\bf FIM}.
Since primitive recursive codes for finite sequences of natural numbers are available in {\bf HA} and hence in {\bf FIM},
to prove the lemma we need only derive a contradiction from
the assumption that $\forall x [\forall y \neg R(x,y) \vee \exists y R(x,y)]$ is derivable in {\bf FIM} from the
universal closures of finitely many instances $\forall x P_i(x,z) \vee \neg \forall x P_i(x,z)$, $1 \leq i \leq k$,
of $\Pi^0_1$-LEM, where $R(x,y)$ is $T(x,x,y)$ and the $P_i(x,z)$ are quantifier-free.  Assume such a derivation
exists, and let $D(z)$ abbreviate $\bigwedge_{i=1}^k (\forall x P_i(x,z) \vee \neg \forall x P_i(x,z))$.
Then by the deduction theorem, {\bf FIM} proves
$${\rm (i)} ~\forall z D(z) \rightarrow \forall x [\forall y \neg R(x,y) \vee \exists y R(x,y)]\enspace.$$
We can construct a purely arithmetical formula $E(w,z)$, with no $\exists$ and no $\vee$, such that {\bf FIM} proves
\begin{enumerate}
\item[]{${\rm (ii)} ~E(w,z) \leftrightarrow \neg \neg E(w,z)~ \mbox{and}$}
\item[]{${\rm (iii)} ~E(\overline{\sigma}({\bf k}),z) \leftrightarrow \\ ~~~~~~\left[\bigwedge_{i=1}^k
( \{\forall x P_i(x,z): \sigma(i\dot{-}1)>0\} \cup \{\neg \forall x P_i(x,z): \sigma(i\dot{-}1)=0\} ) \right]$}
\end{enumerate}
whence
$${\rm (iv)} ~\forall z \left[ D(z) \leftrightarrow
\exists \sigma \in {^{\omega}2}\; E(\overline{\sigma}({\bf k}),z)  \right]~$$
and so
$${\rm (v)} ~\forall z \exists \sigma \in {^{\omega}2}\; E(\overline{\sigma}({\bf k}),z)
\rightarrow \forall x [\forall y \neg R(x,y) \vee \exists y R(x,y)] \enspace.~$$
The countable axiom of choice, which is an axiom schema of {\bf FIM}, gives
$${\rm (vi)} \forall z \exists \sigma \in {^{\omega}2} \; E(\overline{\sigma}({\bf k}),z) \leftrightarrow
\exists \tau \forall z (\lambda t. \tau((z,t)) \in {^{\omega}2} \wedge E(\overline{\lambda t. \tau((z,t))}({\bf k}),z))$$
and hence
$${\rm (vii)} ~\forall \tau \in {^{\omega}2}[\forall z \; E(\overline{\lambda t. \tau((z,t))}({\bf k}),z)
\rightarrow \forall x [\forall y \neg R(x,y) \vee \exists y R(x,y)]]~ $$
where neither $x$ nor $y$ is free in the hypothesis, so also
$${\rm(viii)} ~\forall x \forall \tau \in {^{\omega}2}[\forall z \; E(\overline{\lambda t. \tau((z,t))}({\bf k}),z)
\rightarrow \exists y [\forall y \neg R(x,y) \vee R(x,y)]]~ $$
with a stable hypothesis.  Applying ($\ast$), {\bf FIM} proves
$${\rm(ix)} ~\forall x \forall \tau \in {^{\omega}2} \exists y
[\forall z \; E(\overline{\lambda t. \tau((z,t))}({\bf k}),z) \rightarrow [\forall y \neg R(x,y) \vee R(x,y)]]\enspace.$$
The classically false form of Brouwer's Fan Theorem ($^*$27.7 in \cite{kleenevesley}), followed by the obvious counting
argument, allows us to conclude from (ix) that {\bf FIM} proves
$${\rm(x)} ~\forall x \exists m \forall \tau \in {^{\omega}2} [\forall z \, E(\overline{\lambda t. \tau((z,t))}({\bf k}),z)
 \rightarrow \exists y_{y \leq m} [\forall y \neg R(x,y) \vee R(x,y)]]~$$
and hence
$${\rm(xi)} ~\forall x \exists m [\forall z \exists \sigma \in {^{\omega}2} \; E(\overline{\sigma}({\bf k}),z)
\rightarrow \exists y_{y \leq m} [\forall y \neg R(x,y) \vee R(x,y)]] ~$$
or equivalently
$${\rm (xii)} ~\forall x \exists m \left[\forall z D(z)
\rightarrow \exists y_{y \leq m} [\forall y \neg R(x,y) \vee R(x,y)] \right]\enspace.$$
But then by Kleene's Rule {\bf FIM} proves
$${\rm (xiii)} ~\forall x \exists y \left( T({\bf f},x,y) \wedge \left(\forall z D(z) \rightarrow
 \exists z_{z \leq U(y)} [\forall y \neg T(x,x,y) \vee T(x,x,z)] \right) \right)~ $$
for some natural number $f$, and hence
$${\rm (xiv)} ~\forall z D(z) \rightarrow \exists f F(f)$$
where $F(f)$ is $\forall x \exists y \left( T(f,x,y) \wedge
[\forall z_{z \leq U(y)} \neg T(x,x,z) \rightarrow \forall y \neg T(x,x,y)] \right)$.
Lemma 2 and (xiv) together now imply that {\bf FIM} proves $\neg \forall z D(z)$, which is impossible since {\bf PA}
is consistent with {\bf FIM}.

\vskip 0.1cm

{\em Theorem 1}.  (a)  Each of the arithmetical principles $\Sigma^0_1$-LEM, $\Sigma^0_1$-DNE is independent relative to
the arithmetical principle $\Pi^0_1$-LEM over {\bf FIM}.

(b)  For every $n \geq 1$:  Each of the arithmetical principles $\Sigma^0_{n+1}$-LEM, $\Sigma^0_{n+1}$-DNE is independent
relative to the arithmetical principle $\Pi^0_{n+1}$-LEM over {\bf FIM} + $\Sigma^0_n$-LEM.

\vskip 0.1cm

{\em Proof}.  (a) follows from Lemmas 1-3.  To prove (b) for $n \geq 1$, we need to generalize the lemmas.  Since
$\Pi^0_{n+1}$-LEM implies $\Sigma^0_n$-DNE and $\Sigma^0_n$-LEM, Lemma 1 holds with $\Pi^0_{n+1}$ and $\Sigma^0_{n+1}$
in place of $\Pi^0_1$ and $\Sigma^0_1$ respectively.  Lemma 2 holds with $T^Q$ in place of $T$, where $Q$ is any
$\Sigma^0_n$ predicate.

For Lemma 3 with {\bf FIM} + $\Sigma^0_n$-LEM in place of {\bf FIM}, and $\Pi^0_{n+1}$ and $\Sigma^0_{n+1}$ in place
of $\Pi^0_1$ and $\Sigma^0_1$, take $R(x,y)$ to be the complete predicate for arithmetical $\Pi^0_n$.
Each $P_i(x,z)$ (now $\Sigma^0_n$) is equivalent in {\bf HA} + $\Sigma^0_n$-LEM to its G\"odel-Gentzen
negative translation, so we may use these in defining $E(w,z)$. {\bf FIM} + $\Sigma^0_n$-LEM satisfies $(\ast)$
because $\Sigma^0_n$-LEM proves its own Friedman translation by a stable formula.  The step corresponding to
(xii) $\Rightarrow$ (xiii)
is justified by Theorem 50(b) and Corollary 57 in \cite{redmonograph}, and the contradiction
follows because {\bf PA} is consistent with {\bf FIM} + $\Sigma^0_n$-LEM.

\vskip 0.1cm

{\it Corollary}.  All the derivability and relative independence results over {\bf HA} established by Akama, Berardi,
Hayashi and Kohlenbach among the purely arithmetical principles $\Delta^0_{n+1}$-LEM, $\Pi^0_{n+1}$-LEM,
$\Sigma^0_{n+1}$-DNE and $\Sigma^0_{n+1}$-LEM hold also over {\bf FIM}, for every $n \geq 0$.

\vskip 0.1cm

{\it Proof}. The relative derivability results are preserved because {\bf HA} is a subsystem of {\bf FIM}.
$\Sigma^0_{n+1}$-LLPO is independent relative to $\Sigma^0_{n+1}$-DNE over {\bf FIM} because every theorem of
{\bf FIM} + $\Sigma^0_{n+1}$-DNE is classically realizable by a $\Delta^0_n$ function, while $\Sigma^0_{n+1}$-LLPO is not.
Hence also $\Pi^0_{n+1}$-LEM and $\Sigma^0_{n+1}$-LEM are independent relative to $\Sigma^0_{n+1}$-DNE over {\bf FIM}.

The theorem takes care of the other cases.   For example,
$\Sigma^0_{n+1}$-DNE is independent relative to $\Delta^0_{n+1}$-LEM over {\bf FIM} by the theorem, because
{\bf FIM} + $\Pi^0_{n+1}$-LEM proves $\Delta^0_{n+1}$-LEM but not $\Sigma^0_{n+1}$-DNE.

\vskip 0.1cm

{\em Open Questions}?  I do not know whether $\Pi^0_{n+1}$-LEM is independent relative to $\Sigma^0_{n+1}$-LLPO over
{\bf FIM}.  Lifschitz realizability cannot be used here because {\bf FIM} includes countable and continuous choice
principles.  I also do not know whether $\Delta^0_{n+1}$-LEM is independent of $\Sigma^0_n$-LEM over {\bf FIM}.
Classically, $\Delta^0_1$-LEM is realizable, $_{\sf S}$realizable and $^{\sf G}$realizable so these standard methods do
not give independence even for $n=0$.

\section{How Double Negation Changes the Picture}

Let $\neg \neg \forall$($\Sigma^0_n$-LEM) abbreviate the double negation of the universal closure of arithmetical
$\Sigma^0_n$-LEM, and similarly for the other principles.  For each $n \geq 0$ the weaker principles behave, with
respect to relative independence over {\bf HA}, very much like the stronger ones.

\vskip 0.1cm

{\em Theorem 2}.  Over {\bf HA}, for each $n \geq 1$:

(a) $\neg \neg \forall$($\Sigma^0_n$-LEM) entails $\neg \neg \forall$($\Pi^0_n$-LEM).

(b) $\neg \neg \forall$($\Pi^0_n$-LEM) entails $\neg \neg \forall$($\Delta^0_n$-LEM), but not conversely.

(c) $\neg \neg \forall$($\Sigma^0_n$-LEM) entails $\neg \neg \forall$($\Sigma^0_n$-DNE), but not conversely.

(d) $\neg \neg \forall$($\Sigma^0_n$-DNE) entails $\neg \neg \forall$($\Delta^0_n$-LEM), but not conversely.

(e) $\neg \neg \forall$($\Sigma^0_n$-DNE) does not entail $\neg \neg \forall$($\Pi^0_n$-LEM).

\vskip 0.1cm

{\em Proof}.  Only the relative independence results require comment.  Classical number-realizability relativized to
$\Delta^0_n$ shows that {\bf HA} + $\Delta^0_n$-LEM does not prove $\neg \neg \forall$($\Pi^0_n$-LEM), and that
{\bf HA} + $\Sigma^0_n$-DNE proves neither $\neg \neg \forall$($\Sigma^0_n$-LEM) nor $\neg \neg \forall$($\Pi^0_n$-LEM).
To show {\bf HA} + $\Delta^0_1$-LEM does not prove $\neg \neg \forall$($\Sigma^0_n$-DNE) use modified number-realizability
relativized to $\Delta^0_n$.

\vskip 0.1cm

Does {\bf HA} + $\Pi^0_n$-LEM or {\bf FIM} + $\Pi^0_n$-LEM prove either $\neg \neg \forall$($\Sigma^0_n$-DNE) or
$\neg \neg \forall$($\Sigma^0_n$-LEM)?  I do not know.

Most of Theorem 2 extends to {\bf FIM}, using $^{\Delta^0_n}$realizability (a generalization of the $^{\sf G}$realizability
in \cite{cantherebeno}) for the nonderivabilities in (b) and (d).  However, $\neg \neg \forall$($\Sigma^0_n$-DNE) is
interderivable with $\neg \neg \forall$($\Sigma^0_n$-LEM) over {\bf FIM}, by the following result.

\vskip 0.1cm

{\em Theorem 3}.
(a)  Over {\bf FIM}, and hence over {\bf HA}, each original principle (possibly excepting
$\Sigma^0_1$-DNE and $\Delta^0_1$-LEM) is strictly stronger than its doubly negated closure.

(b)  {\bf FIM} + $\Sigma^0_n$-DNE proves $\neg \neg \forall$($\Sigma^0_n$-LEM), for $n \geq 1$.

(c)  {\bf HA} + $\Sigma^0_n$-DNE does not prove $\neg \neg \forall$($\Sigma^0_n$-LEM).

\vskip 0.1cm

{\em Proofs}.  Each doubly negated closure is classically function-realizable, while $\Sigma^0_1$-DNE and $\Delta^0_1$-LEM
are the only original principles with this property, so (a) holds.

By an argument essentially due to Solovay,
{\bf FIM} + $\Sigma^0_n$-DNE proves $\neg \neg \forall$($\Sigma^0_n$-LEM) for every $n \geq 1$.
The proof in \cite{myhierarchiespaper} using an analytical version of Markov's Principle can be paraphrased
to give the result for the arithmetical principles from arithmetical $\Sigma^0_n$-DNE, so (b) holds also.
Finally, (c) follows from the proof of Theorem 2(c).

\end{document}